\documentclass[aps,prl,singlecolumn,preprintnumbers,amsmath,amssymb,notitlepage]{revtex4-1}
\usepackage{graphicx}
\usepackage{subfigure}
\usepackage{epsfig}
\usepackage{dcolumn}
\usepackage{bm}
\usepackage{ulem}
\usepackage{color}
\usepackage{amsmath}
\usepackage{float}
\usepackage[colorlinks]{hyperref}
\hypersetup{
	citecolor = cyan
}

\def\be{\begin{equation}}       \def\ee{\end{equation}}
\def\bea{\begin{eqnarray}}      \def\eea{\end{eqnarray}}
\def\ba{\begin{array}}
\def\ea{\end{array}}
\def\bnum{\begin{enumerate} }
\def\enum{\end{enumerate}}

\def\=>{\Rightarrow}
\def\>{\rightarrow}

\def\eye2{Fathbb{I}}

\def\Tr{\mathrm{Tr}}

\renewcommand{\>}{\rangle}

\newcommand{\p}{\partial}

\newcommand{\al}[1]{\begin{align}#1\end{align}}
\newcommand{\eq}[2]{
	\begin{equation}
	#1 \label{#2}
	\end{equation}
}

\newcommand{\pf}[2]{\frac{\partial #1}{\partial #2}}

\newcommand{\df}[2]{\frac{\mathrm{d} #1}{\mathrm{d} #2}}

\begin{document}
\title{Automatic Differentiation for Complex Valued SVD}
\author{Zhou-Quan Wan}
\author{Shi-Xin Zhang}
\email{zsx16@mails.tsinghua.edu.cn}
\affiliation{Institute for Advanced Study, Tsinghua University, Beijing 100084, China}

\begin{abstract}
In this note, we report the back propagation formula for complex valued singular value decompositions (SVD).  This formula is an important ingredient for a complete automatic differentiation(AD) infrastructure in terms of complex numbers, and it is also the key to understand and utilize AD in tensor networks.
\end{abstract}

\date{\today}
\maketitle

\section{Introduction}

Automatic differentiation(AD) evaluates derivatives or gradients of any functions specified by computer programs\cite{Biggs2000}. It is implemented by propagating derivatives of primitive operations via chain rules. Such approach is different from classical symbolic or numerical differentiations.
Symbolic differentiation faces the challenge of converting a complicated computer program into expressions, while numerical differentiation faces the difficulty of numerical errors in the discretization. 
Besides, both symbolic and numerical methods have problems in calculating higher order derivatives and are also slow at computing gradients with respect to lots of inputs variables, e.g. in the case for gradient-based optimization algorithms. 
AD solves all of these drawbacks and is emerging as a new programming paradigm which is extensively utilized in machine learning(ML) community and serves as the most important infrastructure for ML libraries.

AD find applications in various other fields as well. Specifically, AD has been applied to computational physics where the interplay between tensor networks and AD has been studied recently\cite{Liao2019, Torlai2019}. In the context of quantum physics, complex number plays an important role which is often overlooked by ML community. Therefore, although there are many works on AD for real valued linear algebra\cite{Giles2008, Townsend2016, Seeger2017}, there are very few reports on AD for complex-valued linear algebra\cite{Hubig2019, Liu2019}. And as mentioned by both \cite{Torlai2019} and \cite{Hubig2019}, reverse mode AD formula for complex valued SVD is still missing. AD formula for complex valued SVD would allow us to directly utilize AD in tensor network context where tensors are complex valued and operations include SVD. AD knowledge on complex valued SVD is necessary to explore AD on tensor networks for complicated models as well as other scenarios where SVD lies within the forward functions for numerical stability or information extraction.

\section{Results}

In this note, we give the formula and derivation of reverse mode AD for complex valued SVD.

For square matrix SVD $A=USV^\dagger$, with the final loss function $L$ being real, the back propagation formula is given by the following:
\begin{align}
	\bar{\mathcal {A} }= U\bar{\mathcal{S}}V^\dagger+U(\mathcal{J}+\mathcal{J}^\dagger)SV^\dagger + US(\mathcal{K}+\mathcal{K}^\dagger) V^\dagger+\frac{1}{2}US^{-1}(\mathcal{L}^\dagger-\mathcal{L})V^\dagger,
	\label{square}
\end{align}
where $\mathcal{J}=F\circ (U^\dagger \bar{\mathcal{U}}), \mathcal{K}=F\circ (V^\dagger \bar{\mathcal{V}}), \mathcal{L}=I\circ (V^\dagger \bar{\mathcal{V})}$ and
$$
F_{ij}=\left\{ \begin{aligned}\frac{1}{s_j^2-s_i^2}~(i\neq j)\\0~~~~~~ (i=j)\end{aligned}\right . 
$$

The main contribution of this work is the fourth term in \eqref{square} as
\eq{\bar{\mathcal{A}}_d=\frac{1}{2}US^{-1}(\mathcal{L}^\dagger-\mathcal{L})V^\dagger.}{diagonal}

This term is new and has no counterpart in real valued SVD.

\section{Derivation}

{\bf Basics on complex calculus and some notations:}

The partial derivatives for function $f(z)$ where $z=x+iy$ is defined as:
\eq{\pf{f}{z}=\frac{1}{2}(\pf{f}{x}-i \pf{f}{y}), ~~~ \pf{f}{z^*}=\frac{1}{2}(\pf{f}{x}+i \pf{f}{y}).}{derivative}

Such form of derivative is justified as 
\eq{df=\pf{f}{z}dz+\pf{f}{z^*}dz^*.}{diff}

If $f$ is a holomorphic function of $z$, then $\pf{f}{z^*}=0$. Only for these holomorphic $f$, we have well-defined derivatives for complex function as 
\eq{\df{f}{z}=\pf{f}{z}.}{arg2}

In this note, we only focus on real valued functions depending on complex variables. This is the case as we can only optimize a real value; optimization problem makes no sense for complex valued object. Besides, for physics settings, all observables give real valued results. Note real valued $f$ cannot be holomorphic unless $f$ is a constant independent of $z$.

For optimization problems and gradient descent approach, what we really care about is gradients instead of derivatives. For real variable $x$ and real function $f(x)$, it is coincident that $\nabla_x f = \pf{f}{x}$. However this relation doesn't hold for complex variable functions. Instead, for complex variable functions, the gradient is defined as
\eq{\nabla_z f(z)= 2\pf{f}{z^*}=\pf{f}{x}+i\pf{f}{y}.}{grad}
If $f$ is real valued function, we further have $\nabla_z f = 2(\pf{f}{z})^*$. See \cite{Hunger2007} for more details on the calculus of complex and matrix based functions.

 In the following, we only focus on reverse mode AD only. We use $L$ to represent the final real valued loss function, and use the following notation for derivatives and gradients:
 \eq{\bar{A}=\pf{L}{A}~~~\bar{\mathcal{A}}=2\pf{L}{A^*}=\nabla_AL.}{notat} 
 
 $\circ$ stands for elements multiplication between two matrices. The following facts will be utilized in the proof below:
 \eq{\Tr[A(C\circ B)]=\Tr[(C^T\circ A)B],}{arg2}
 \eq{(A\circ B)^T=A^T\circ B^T.}{arg2}
 
 $I$ denotes the identity matrix while $\bar{I}$ represents a square matrix with diagonal zero and all off-diagonal elements 1.
 
{\bf Parallel proof with real case:}

For the first three terms in \eqref{square},  we follow the similar proof for real valued SVD with special attention on unique features for complex functions. The SVD is defined as:
\eq{A=USV^\dagger,}{svd}
where A is the input matrix while U, S, V is the output matrix. S only has positive real elements in the main diagonal, and U, V are unitary matrix: $U^\dagger U=I, V^\dagger V=I$.

A direct differentiation on \eqref{svd} gives:
\eq{dA=dUSV^\dagger +UdSV^\dagger+USdV^\dagger.}{dsvd}
Let $dC=U^\dagger dU, dD=V^\dagger dV$, according to the unitary property of U and V, we have 
\eq{dC=-dC^\dagger ~~ dD=-dD^\dagger.}{antisymmetry}
\eqref{antisymmetry} indicates that the diagonal elements in dC and dD are {\bf pure imiginary}. (not zero!) Plug dC and dD into \eqref{dsvd}, we have
\eq{dP=dCS+dS-SdD,}{main}
where $dP=U^\dagger dAV$.

We have several relations from real and imaginary diagonal part as well as off-diagonal parts from \eqref{main}, respectively. (Note that S and dS is real diagonal matrix and dC, dD have pure imaginary diagonal elements.)
\eq{dS=I\circ (\frac{dP+dP^\dagger}{2}),}{dS}
\eq{I\circ(dC-dD)=(I\circ (\frac{dP-dP^\dagger}{2}))S^{-1},}{dCd}
\eq{\bar{I}\circ(dP S+S dP^\dagger)=dC S^2 - S^2 dC,}{dCo}
\eq{\bar{I}\circ(S dP+dP^\dagger S)=dD S^2 - S^2 dD.}{dDo}

Our proof is based on the following relation given by chain rules:
\eq{dL=\Tr[\bar{U}^TdU+\bar{V}^TdV+\bar{S}^TdS+c.c.]=\Tr[\bar{A}^TdA +c.c.]}{dL} 
and our aim is to find the analytical formula for $\bar{A}$ in terms of U, S, V and $\bar{U}, \bar{S}, \bar{V}$. $dL$ in \eqref{dL} can be decomposed as
\begin{align}
	dL&=\Tr[\bar{S}^TdS+\bar{U}^TU(I\circ dC)+\bar{U}^TU(\bar{I}\circ dC)+\bar{V}^TV(I\circ dD)+\bar{V}^TV(\bar{I}\circ dD)+c.c.]\\
	&=\Tr[\bar{A}^T_{s}dA+\bar{A}^T_{ud}dA+\bar{A}^T_{uo}dA+\bar{A}^T_{vd}dA+\bar{A}^T_{vo}dA+c.c.].
\end{align}

For $\bar{A}_s$ part, we have:
\eq{\Tr[\bar{S}^T (I\circ dP)+c.c.]=\Tr[\bar{S}^T U^\dagger dA V+c.c.]=\Tr[V\bar{S}^TU^\dagger dA+c.c.].}{arg2}
Namley, 
\eq{\bar{A}_s=U^*\bar{S}V^T.}{As}

For $\bar{A}_{uo}$ part, note $\bar{I}\circ dC=F\circ (dP S+SdP^\dagger)$, we have:
\al{&\Tr[\bar{U}^TU(\bar{I}\circ dC)+c.c.]=\Tr[\bar{U}^T U F\circ (dP S+ S dP^\dagger)+c.c.]=\\
&	=\Tr[(-F\circ \bar{U}^TU)(U^\dagger dA VS +SV^\dagger dA^\dagger U)+c.c.]=\Tr[VSJ^TU^\dagger dA+V SJ^* U^\dagger dA+c.c.]
	,}
where $J=F\circ (U^T\bar{U})$. From this, we conclude that
\eq{\bar{A}_{uo}=U^* (J+J^\dagger)SV^T.}{Auo}

Similarly, we have:
\eq{\bar{A}_{vo}=U^* S (K+K^\dagger)V^T,}{Avo}
where $K=F\circ (V^T\bar{V})$.

The three terms \eqref{As}, \eqref{Auo} and \eqref{Avo} as well as the proof above are very similar with real valued SVD. The difference here is that dC and dD have pure imaginary diagonals instead of zero diagonals. Consequently, there are two more differentiation terms for $dA$, corresponding to $dA_{ud}$ and $dA_{vd}$, originating from the differentiation of diagonal part in dC and dD. The two terms are totally new in complex context, and we will explore how to compute them in the following sections.

{\bf Gauge freedom:}

There are some extra freedom for the choice of U and V in SVD as we can see from the following relation:
\eq{A=USV^\dagger=U\Lambda \Lambda^\dagger S\Lambda \Lambda^\dagger V^\dagger=U'SV'^\dagger ,}{gauge}
where $U'=U\Lambda, V'=V\Lambda$, $\Lambda$ is a diagonal matrix with diagonal elements in the form of a phase  $e^{i\theta}$. Therefore, the loss function $L$ must be gauge invariant, i.e. $L(U, S, V)=L(U\Lambda, S, V\Lambda)$. 

Now considering gauge transformation $\Lambda$ on $dL$:
\al{dL(U, S, V)=dL(U\Lambda, S, V\Lambda)=\Tr[\bar{S}^TdS+\bar{U}^TdU+\bar{V}^TdV+c.c.]\nonumber\\
=\Tr[\bar{S}^TdS+\overline{U\Lambda}^T d(U\Lambda)+\overline{V\Lambda}^Td(V\Lambda)+c.c.].}
Gauge invariance indicates that:
\eq{\Tr[\bar{\Lambda}^Td\Lambda+\bar{\Lambda}^\dagger d\Lambda^*]=\Tr[(\overline{U\Lambda}^TU+\overline{V\Lambda}^TV)d\Lambda+(\overline{U\Lambda}^\dagger U^*+\overline{V\Lambda}^\dagger V^*)d\Lambda^*]=0.}{invp}
Since $\Lambda$ is diagonal matrix with phase elements $e^{i\theta_j}, (j=1,2..n)$, the differential of it around $\Lambda=I$ gives that $d\Lambda^* = -d\Lambda$. Namely $d\Lambda$ is an imaginary diagonal matrix with each diagonal element differential $d\theta_j$ independent on each other and can vary freely.

Taking $\Lambda=I$  we have, 
\eq{\Tr[(\bar{U}^TU+\bar{V}^TV-c.c.)d\Lambda]=0.}{}
Based on the fact that each diagonal term of $d\Lambda$ is arbitrary, 
\eq{(\bar{U}^TU+\bar{V}^TV-c.c.)\circ I=0.}{inv}
We have obtained the important equation \eqref{inv} from gauge invariance.

It should be emphasized that gauge invariance objective $L$ doesn't directly imply $\bar{\Lambda}=\frac{\p L}{\p \Lambda}=0$.
And $\bar{\Lambda}$ can be indeed nonzero even if $L$ is gauge invariant. For example, $L=u^*_{00}u_{00}$ is invariant under phase gauge. For this objective, $\bar{V}=0$ and $\bar{U}$ also has all zero elements except $\bar{u}_{00}$. It is obvious $\bar{\Lambda}_{0i}=[\bar{U}^TU+\bar{V}^TV]_{0i}=\bar{u}_{00}u_{0i}$ is not zero in general.

{\bf New differentiation term in complex case:} 

Now considering the remaining term $\bar{A}_{ud}+\bar{A}_{vd}$, we have:
\eq{\Tr[\bar{U}^TU(I\circ dC)+\bar{V}^TV(I\circ dD)+c.c.]=\Tr[(\bar{U}^TU+\bar{V}^TV)(I\circ dC)+c.c.]+\Tr[\bar{V}^TV(I\circ (dD -dC))+c.c.].}{dec}
The first term on rhs is zero due to gauge invariance \eqref{inv} by noting $I\circ dC$ is pure imaginary,
\eq{
\begin{split}
\Tr[(\bar{U}^TU+\bar{V}^TV)(I\circ dC)+c.c.]&=\Tr[(\bar{U}^TU+\bar{V}^TV-c.c.)(I\circ dC)]\\
&=\Tr[(\bar{U}^TU+\bar{V}^TV-c.c.)\circ I) dC]=0.
\end{split}
}{}
The remaining one can be replaced by \eqref{dCd}:
\eq{\Tr[(\bar{A}_{ud}+\bar{A}_{vd})^TdA+c.c.]=\frac{1}{2}\Tr[\bar{V}^TV(I\circ (dP^\dagger-dP)S^{-1})+c.c.]
=	\frac{1}{2}\Tr[V (L^* S^{-1}-S^{-1}L^T)U^\dagger dA+c.c.],}{arg2}
where $L=I\circ (V^T\bar{V})$.

From this, we have:
\eq{\bar{A}_d\equiv \bar{A}_{ud}+\bar{A}_{vd}=\frac{1}{2} U^*  (L^\dagger S^{-1} -S^{-1}L) V^T.}{Aud} 

Now we have derived all contributions to $\bar{A}$, the sum of \eqref{As}, \eqref{Auo}, \eqref{Avo} and \eqref{Aud} gives the back propagation formula of derivatives for complex valued SVD. 

It is worth noting the formula form of \eqref{Aud} and thus the formula form of reverse mode AD for complex valued SVD is not unique due to gauge freedom in SVD since we can similarly keep $\Tr[\bar{U}^TU(I\circ (dC-dD))]$ instead of $\Tr[\bar{V}^TV(I\circ (dD-dC))]$ in the calculation. But different formulas give the same numerical value due to gauge invariance \eqref{inv}.

In practice, our aim is to obtain gradients of $L$, since $\bar{\mathcal{A}}=2\bar{A}^*$ for real $L$, we can easily transform the relation for derivatives to the one for gradients, and the final result is \eqref{square} \cite{Wan2019}. Note how $\mathcal{J}, \mathcal{K}, \mathcal{L}$ changed accordingly from $J, K, L$.

\section{Discussions}

The results can directly apply to SVD on rectangular matrix A, with two extra terms similar to the real SVD case, namely

\eq{\bar{\mathcal{A}}=\bar{\mathcal{A}}_{square}+(I-UU^\dagger)\bar{\mathcal {U}}S^{-1}V^\dagger +US^{-1}\bar{\mathcal {V}}^\dagger (1-VV^\dagger), }{rec}
where $\bar{ \mathcal{A}}_{square}$ is defined as \eqref{square}. In other words, for general complex valued SVD, we only need to add one extra term beyond real SVD case, which is just $\bar{\mathcal{A}}_d$ in \eqref{diagonal}.

It is worth noting that if the loss function has no dependence on U or V, i.e. $L=L(U,S)$ or $L=L(V,S)$, the contribution $\bar{A}_{d}=0$ following \eqref{inv}. In other words, to test the correctness of backprop for complex valued SVD, the loss function must depend on both U and V at the same time. Otherwise, the original formula for SVD without $\bar{A}_{d}$ still works. Therefore, one must devise a test loss function $L$, which is real valued, gauge invariant, and cannot decouple to two gauge invariant parts depending on U, S and S, V only. An example loss function might be something like $\mathcal{R}(U_{00}V_{00}^*)$.

%It is also possible to derive the forward mode AD formula for complex valued SVD operation. The simplest way to do this is by gauge fixing, i.e. by taking $I\circ dC =0$. We can then obtain forward AD formula from \eqref{dS} to \eqref{dDo}. Different choices of gauge fixing give different forward AD formulas and $U', S', V'$ can even be different in numerical values. This is on the contrary to reverse mode AD, where different reverse mode AD formulas are possible but they all give the same numerical value for $\bar{A}$. But forward mode AD is still reasonable since $L'=\frac{\p L}{\p U}U'+\frac{\p L}{\p S}S'+\frac{\p L}{\p V}V'$ is the same by gauge invariance of $L$, though $U', V', S'$ can be different.

AD formula for complex valued SVD has great potential for applications in various fields ranged from machine learning to computational physics. For example, AD for complex valued SVD may advance the development of relevant optimization algorithms on tensor networks\cite{Liao2019}.
\newline
\newline
{\it Acknowledgement}: 
We thank Jin-Guo Liu for helpful discussions on this topic.

\end{document}